\input epsf
\epsfverbosetrue
\input amssym.def
\magnification=1200

\baselineskip=14pt

\def\qed{{$\vrule height4pt depth0pt width4pt$}}
\def\ss{{\smallskip}}
\def\ms{{\medskip}}
\def\bs{{\bigskip}}
\def\ni{{\noindent}}

\def\t{{\tau}}
\def\p{{\partial}}

\def\g{{\gamma}}

\def\e{{\epsilon}}
\def\l{{\lambda}}
\def\b{{\beta}}
\def\Z{\Bbb Z}

\def\<{{\langle}}
\def\>{{\rangle}}

\def\D{{\Delta}}

\def\<{{\langle}}
\def\>{{\rangle}}

\def\m1{{\quad ({\rm mod\ 1})}}
\def\L{{\Lambda}}

\def\={\buildrel \cdot \over =}


\centerline{\bf $3$-MANIFOLDS, TANGLES AND PERSISTENT INVARIANTS}

\bs
 
\centerline{J\'OZEF H. PRZYTYCKI}
\centerline{Department of Mathematics, The George Washington
University}
\centerline{Washington, DC 20052, USA}
\centerline{e-mail: przytyck@gwu.edu}\ms

\centerline{DANIEL S. SILVER and SUSAN G. WILLIAMS}
\centerline{Department of Mathematics and Statistics, University
of South Alabama}
\centerline{Mobile, AL  36688-0002, USA}
\centerline{e-mail: silver@jaguar1.usouthal.edu,
swilliam@jaguar1.usouthal.edu} \bs

\footnote{} {Second and third authors partially supported by NSF
grant DMS-0071004.}
\footnote{}{2000 {\it Mathematics Subject Classification.}  
Primary 57M25; secondary 11R06.}

\noindent {\narrower\narrower\smallskip\noindent  
ABSTRACT: Given a compact, connected,
oriented $3$-manifold $M$ with boundary,
and epimorphism
$\chi$ from $H_1M$ to a free abelian group $\Pi$,
two invariants
$\b$, $\t \in \Z\Pi$ are defined.  If $M$
embeds in another such $3$-manifold $N$ such that 
$\chi_N$ factors through $\chi$, then the product $\b\t$ divides
$\D_0(H_1\tilde N)$.  

A theorem of D. Krebes concerning 
$4$-tangles embedded in links arises as a special 
case. Algebraic and skein theoretic generalizations 
for $2n$-tangles provide invariants that persist in the
corresponding invariants of links in which they embed. An
example is given of a virtual $4$-tangle for which
Krebes's theorem does not hold.


\smallskip}
\ms


\ni {\bf 1. Introduction.}\bs

Suppose that a $2n$-tangle $t$ embeds in a link $\ell$.
It is natural to ask which invariants of $t$ necessarily persist as invariants of
$\ell$. In [{\bf Kr99}] D. Krebes considered the case that $t$ is a
$4$-tangle. There he proved that any positive integer dividing the determinants of both the numerator closure and the denominator closure of
the tangle also divides the determinant of $\ell$. 

Krebes's argument is a blend of combinatorics and topology. In
[{\bf Ru00}] D. Ruberman gave another proof of Krebes's
theorem using purely topological techniques. Ruberman exploited
a well-known relationship between the determinant of a link and 
the first homology of its $2$-fold cyclic branched cover. From his
perspective Krebes's theorem is a result about invariants of 
compact, oriented $3$-manifolds that persist as invariants of
rational homology $3$-spheres in which they embed. 

We extend Ruberman's methods in order to prove a generalization of
Krebes's theorem for $2n$-tangles (see Theorem 2.5). Given a
$2n$-tangle $t$ we define two invariants $\t, \b \in 
\Z\Pi$, where $\Pi$ is a free abelian group. The 
rank of $\Pi$, denoted by $d$,  depends on the category in which we work: If we
color and orient the components, then $d$ can be chosen to be 
the number of components of $t$. Then whenever $t$
embeds in an oriented link $\ell$ such that distinct components of
$t$ lie in different components of $\ell$, the product $\t\b$
divides the multivariable Alexander polynomial of $\ell$. At the
other extreme, we may choose to ignore both the colors and
orientations of $t$. In that case $d=0$ and $\t, \b$ are integers. If
$t$ embeds in an unoriented  link $\ell$, then the product $\t\b$
divides the determinant of $\ell$. The latter statment is seen to be
Krebes's theorem for $4$-tangles.

The second and third authors gave a short, elementary proof of
Krebes's theorem when the divisor is prime. The proof
immediately extends when the divisor is square-free. The argument, based
on the combinatorial technique of Fox coloring, holds in the larger
category of virtual tangles and links [{\bf SW99}]. In Section 4
we consider the virtual category. We give an
example that shows that Krebes's theorem is not valid in the
larger category unless the divisor is square-free. 

Another generalization of Krebes's theorem for $2n$-tangles,
proved using Kauffman bracket skein theory and Temperley-Lieb algebra, is in
[{\bf KSW00}]. We extend this approach to other skein theories in the last
section. 

We are grateful to J. Scott Carter, Mietek D\c abkowski and Seiichi Kamada for
stimulating and helpful discussions. 
\bs

\ni {\bf 2. Persistent invariants of submanifolds.} 
Let $\Pi$ be a free abelian multiplicative group on $d\ge 0$
generators $x_i$. The group ring $\L= \Z\Pi$ is  a Noetherian
unique factorization domain with automorphism $r \mapsto \bar r$
extending the assignment $x_i \mapsto x_i^{-1}$, for all $i$. 

Let $H$ be a finitely generated $\L$-module. The {\it rank}
of $H$ is the dimension of the vector space $Q(\L)\otimes_\L H$, 
where $Q(\L)$ is the field of fractions of $\L$. The {\it
$\L$-torsion submodule} of $H$ is $TH=\{a\in H\mid ra=0
\ {\rm for\ some\ nonzero\ }r\in \L\}$. We denote the {\it Betti module}
$H/TH$ by $BH$. Assume that we have a presentation 
$$\L^p\ {\buildrel R\over \longrightarrow}\ \L^q \to H\to 0.$$
By adding trivial relators, if necessary, we can assume that $q 
\le p$. 
For each $0\le k<q$, the $k${\it th\ elementary\ divisor}
$\D_k(H)$ is  the greatest common divisor of the $(q-k)\times
(q-k)$ subdeterminants of the matrix representing $R$. It is
well defined up to multiplication by a unit in $\L$. By convention,
$\D_k(H)=0$ if $k$ is negative,  while $\D_k(H)=1$ if $k\ge q$. For each $k$, the  polynomial 
$\D_k(H)$ is an invariant of $H$; in particular, it does not depend
on the particular choice of matrix representing $R$. 

Consider a compact $3$-manifold $X$ with
boundary $\p X$ decomposed as the union of two surfaces $\p_+ X$ and $ \p_-X$; if both are nonempty then their intersection should be a $1$-manifold.   Let $\chi:H_1X
\to
\Pi$ be an epimorphism. The map $\chi$ determines an abelian
cover $p:\tilde X\to X$ with deck transformation group
$\Pi$. We denote the preimage $p^{-1}(\p_\pm X)$ by 
$\p_\pm \tilde  X$. The homology
groups $H_*(\tilde X), H_*(\tilde X, \partial \tilde X)$ (integer
coefficients understood) are in fact finitely generated $\L$-modules.
We consider the
composite homomorphism
$$\nabla :H_1\p_+\tilde X\ {\buildrel i_{1+} \over \to}\
H_1\tilde X\  {\buildrel \pi \over
\to}\ BH_1\tilde X,\eqno(2.1)$$ where $i_{1+}$ is the map induced by
inclusion and
$\pi$ is the natural quotient map.   \bs

\ni {\bf Definition 2.1.}
The {\it boundary invariant}\ $\b(X, \p_+X)$ is
$\D_0(BH_1\tilde X/{\rm im}\ \nabla)$. The {\it torsion
invariant}\ $\t(X)$ is
$\D_0(TH_1\tilde X)$.\bs

\ni Although $\t(X)$ and $\b(X,\p_+X)$ depend on $\chi$, 
we omit specific mention of $\chi$ in our notation for the sake of convenience.

Let $M, N$ be compact, connected oriented $3$-manifolds 
with $M\subset N$. Regard $N$ as the union
of $M$ and another compact, oriented $3$-manifold $M'$
with  $M\cap M' = \p_+M$. We assume that 
$\chi$ extends over $H_1N$. The preimage $p^{-1}(M)$ is
connected and can be identified with $\tilde M$. 
Assume that if $\partial \tilde N$ is
nonempty, then each component is noncompact; the assumption
is equivalent to the statement that each component of $\p N$
contains a cycle $z$ such that $\chi([z]) \ne 0$.
 \bs

\ni {\bf Theorem 2.2.} Under the above hypotheses, 
$\b(M, \p_+M) \t(M)$ divides $\D_0(H_1\tilde N)$.

\bs
The proof of Theorem 2.2 when $d>0$ (that is, when the covers are nontrivial) uses Blanchfield duality, which we
review. Let $X$ be a compact, connected $n$-manifold
with boundary $X=\p_+X \cup \p_-X$. If $p: \tilde X
\to X$ is any connected cover of $X$
with covering group
$\Pi$, there there are nondegenerate $\L$-sesquilinear forms
$$\eqalign{BH_p(\tilde X, \p_+\tilde X) \times BH_{n-p}(\tilde X,
\p_-
\tilde X) &\to \L;\cr T_DH_p(\tilde X, \p_+\tilde X) \times
T_DH_{n-p-1}(\tilde X,
\p_- \tilde X) &\to Q(\L)/\L.}$$
Here $T_DH$ denotes the quotient $TH/DH$, where $DH=\{a\in H\mid
r_1a = \cdots = r_qa=0,\ {\rm for\ some\ coprime\ } r_1, \ldots,
r_q \in \L\ (q\ge 2)\}$. Details can be found in [{\bf Bl57}] or
[{\bf Kw96}]. \bs

\ni {\bf Lemma 2.3.} (1) If $0\to A \to B\to C\to 0$ is a short
exact sequence of finitely generated $\L$-modules, then
$\D_0(B) \= \D_0(A)\D_0(C)$, where $\=$ denotes equality in
$\L$ up to multiplication by a unit.  

(2) Let $H$ be a finitely generated $\L$-module. Then $H$ is 
a $\L$-torsion module if and only if $\D_0(H)\ne 0$. More generally,
if $H$ is any finitely generated $\L$-module of rank $r$, then 
$$\D_k(H) \= \cases {$0$ & for $k< r$,\cr \D_{k-r}(TH) & for $k\ge r$}.$$

(3) Let $H$ be a finitely generated $\L$-module. If $D_0$ is a
submodule of $DH$, then $\D_k(H/D_0)= \D_k(H)$ for any $k\ge 0$.

(4) If $0\to A\ {\buildrel f \over \longrightarrow}\ B\ \ 
{\buildrel g \over \longrightarrow}\ C\ \to 0$ is a short exact
sequence of modules over any ring, then for any submodule
$D\subset B$, the following sequence is also exact.
$$0\to A/f^{-1}(D) \ {\buildrel \bar f \over \longrightarrow}\
B/D\
\  {\buildrel \bar g \over \longrightarrow}\ C/g(D) \ \to 0.$$\bs

\ni {\bf Proof.} Lemma 2.3 (1) is well known. A proof can be found
on page 92 of [{\bf Kw96}], for example. The second  and
third statements of Lemma 2.3 are  proved in [{\bf Bl57}] (see
Lemmas 4.3 and 4.10 for (2); for (3) see the proof of Theorem 4.7).
The proof of statement (4) is routine and is left to the reader.
\qed \bs

\ni {\bf Proof of Theorem 2.2.} If $\D_0(H_1\tilde N)=0$,
then the proof of (1) is trivial. Therefore we assume that 
$\D_0(H_1\tilde N)\ne 0$. By Lemma 2.3(2) $BH_1\tilde N =0$.
It follows by Blanchfield duality that $BH_2(\tilde N, \partial
\tilde N)=0$. Thus $H_2(\tilde N, \partial \tilde N)$ is a 
$\L$-torsion module. 

By hypothesis, each component of $\partial \tilde N$ is
noncompact. Hence $H_2\partial \tilde N =0$. From the exact
sequence of the pair $\partial \tilde N \subset \tilde N$:
$$\cdots \to H_2\partial \tilde N \to H_2\tilde N \to
H_2(\tilde N, \partial \tilde N) \to \cdots$$
we see that $H_2\tilde N$ is also a $\L$-torsion module. 

By excision  $H_*(\tilde N, \tilde M)\cong H_*(\tilde M', \p_+ \tilde M)$. Hence
Blanchfield duality pairs $T_DH_2(\tilde N, \tilde M)$ with 
$T_DH_0\tilde(\tilde M', \partial_- \tilde M)$. Since $H_0(\tilde M', \partial_- \tilde M)$ is free, the module $T_DH_0\tilde(\tilde M', \partial_- \tilde M)$ is trivial. Hence $TH_2(\tilde N, \tilde M) \cong DH_2(\tilde N, \tilde M)$. 

Consider now the exact sequence of the pair $\tilde M \subset
\tilde N$:

$$\cdots \to TH_2\tilde N\ {\buildrel k_2 \over \longrightarrow}
\ H_2(\tilde N, \tilde M)\ {\buildrel \partial_2 \over
\longrightarrow}\ H_1\tilde M\ {\buildrel j_1 \over
\longrightarrow}\ H_1\tilde N\ \to\ \cdots$$
From this, $H_1\tilde M/{\rm im}\ \partial_2$ is isomorphic to a
submodule of $H_1\tilde N$. Lemma 2.3 (1) implies that 
$\D_0(H_1\tilde M/{\rm} {\rm im}\ \partial_2)$ divides $\D_0(H_1\tilde N)$. 
It also follows that ${\rm im}\ \partial_2 \cap TH_1\tilde M \subset
DH_1\tilde M$. The reason is the following. If $\partial_2a \in
TH_1\tilde M$, then $ra \in {\rm ker}\ \partial_2$, for some $0\ne
r\in \L$. The exact sequence shows that $ra \in {\rm im}\ k_2$ is a 
torsion element. Consequently $a$ itself is torsion. Hence 
$a \in TH_2(\tilde N, \tilde M) \cong DH_2H(\tilde N, \tilde
M)$, and so $\partial_2 a \in DH_1\tilde M$. 

Consider the canonical short exact sequence:
$$0 \to\ TH_1\tilde M\ \to\ H_1\tilde M\ {\buildrel \pi\over
\longrightarrow}\ BH_1\tilde M\ \to\ 0, $$
where the first map is inclusion and the second is the natural 
quotient projection. By Lemma 2.3 (4):

$$0 \to\ {TH_1\tilde M\over {{\rm im}\ \partial_2 \cap TH_1\tilde M}}\
\to\ {H_1\tilde M\over {\rm im}\ \partial_2}\ {\buildrel
\pi\over
\longrightarrow}\ {BH_1\tilde M\over {\rm im}\ \pi\partial_2}\ \to\
0.$$
By Lemma 2.3 (1) we have
$$\D_0\Bigr({H_1\tilde M\over {\rm im}\ \partial_2}\Bigr)\=
\D_0\Bigr({TH_1\tilde M\over {{\rm im}\ \partial_2 \cap TH_1\tilde
M}}\Bigr)\D_0\Bigr({BH_1\tilde M\over {\rm im}\ \pi\partial_2}\Bigr).$$
Since ${\rm im}\ \partial_2\cap TH_1\tilde M \subset DH_1\tilde M$, by Lemma 2.3 (3) 
$$\D_0\Bigr({TH_1\tilde M\over {{\rm im}\ \partial_2 \cap TH_1\tilde
M}}\Bigr) \= \D_0(TH_1\tilde M),$$
which is the torsion invariant $\t(M)$. 
Also, $BH_1\tilde M/{\rm im} \ \pi\  i_{1+}$ is a quotient of $BH_1\tilde
M/{\rm im}\ \pi\ \partial_2$, and hence the boundary invariant
$$\b(M, \p_+M)
=\D_0\Bigr({BH_1\tilde M\over {\rm im}\ \pi i_1}\Bigr)\ {\rm
divides}\ \D_0\Bigr({BH_1\tilde M\over im\
\pi\partial_2}\Bigr),$$
again using Lemma 2.3 (1). Hence $\b(M,\p_+M)\t(M)$
divides $\D_0(H_1\tilde M/{\rm im}\ \pi \partial_2)$, which we have
previously seen divides $\D_0(H_1\tilde N)$. \qed

\bs
\ni {\bf 3. Application to tangles.} A $2n$-{\it tangle},
for $n$ a positive integer, consists of $n$ disjoint arcs and any
finite number of simple closed curves properly embedded in the
$3$-ball. Two $2n$-tangles are regarded as the same if one can be
transformed into the other by an ambient isotopy of the $3$-ball that
fixes each point on the boundary. As usual we represent
$2n$-tangles by diagrams. Two diagrams represent the same
$2n$-tangle if one can be transformed into the other by a finite
sequence of Reidemeister moves. \bs

\ni {\bf Definition 3.1.} A $2n$-tangle $t$ {\it embeds} in a link
$\ell$ if some diagram for $t$ extends to a diagram for $\ell$. \bs

A $4$-tangle is called simply a {\it tangle}. By joining the top
ends and then the bottom ends one obtains a link $n(t)$,  the 
{\it numerator} of $t$. Joining the left-hand ends and then
right-hand ends produces the {\it denominator} closure $d(t)$.
See Figure 1. \bs

\ni The {\it determinant} ${\rm det}(\ell)$ can be defined in many
ways. It is the absolute value of the one-variable Alexander polynomial
(see below) of $\ell$ evaluated at $-1$. It is also the order of the first homology of the $2$-fold cover of $S^3$ branched over
$\ell$, provided that the group is finite; if not, then the determinant is zero. See for example [{\bf Li97}].
\bs

\epsfxsize=1.5truein
\centerline{\epsfbox{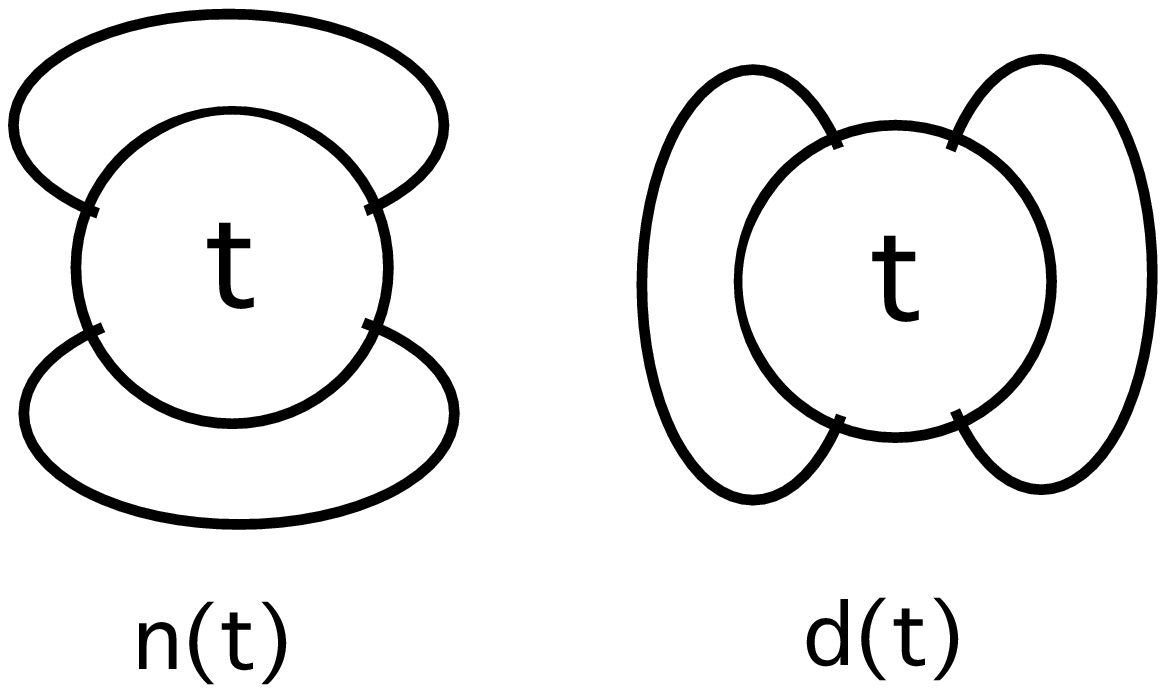}} \bs

\centerline{Figure 1. Diagrams of $n(t)$ and $d(t)$}\bs

\ni {\bf Theorem 3.2.} [{\bf Kr99}] If a tangle $t$ embeds in a
link $\ell$, then the greatest common divisor of ${\rm det}(n(t))$
and ${\rm det}(d(t))$ divides ${\rm det}(l)$. \bs

Krebes's theorem generalizes in various ways. In order to state some
of these, we need some terminology.
If a $2n$-tangle or link has a specified direction for each
component, then it is {\it oriented}. If it has a specified color
for each component, then it is {\it colored}. 

Let $t$ be a colored, oriented $2n$-tangle with exterior
$E_t = B^3-int\ N(t)$. Here $N(t)$ denotes a tubular
neighborhood of $t$. The homology group
$H_1 E_t$ is freely generated by $d$ oriented meridians,
where $d$ is the number of connected components of $t$.

We use Defintion 2.1 to associate invariants $\b, \t$ to the
$3$-manifold $M=E_t$. There are a variety of choices for
$\chi$. If we intend to keep track of both orientations and
colors of $t$, then we consider the isomorphism $\chi_M:
H_1 M \to \Pi \cong \<x_1, \ldots, x_d\mid \ \>$ that maps
the class of the $i$th oriented meridian to  $x_i$. In
this case,
$\tilde M$ is the universal abelian cover of $M$.
Alternatively, we can keep track of orientations but ignore
colors. Then we consider $\chi_M: H_1M \to \Pi\cong
\<x\mid \ \>$, mapping the class of each oriented meridian to $x$. The
covering space $\tilde M$ is sometimes called the ``total
linking number cover.'' 

Krebes's theorem will arise in another way, letting
$M$ be the $2$-fold cyclic cover of $B^3$  branched over $t$ and letting $\Pi$ be the  trivial group. In this case,  $\tilde M$ is equal
to $M$. 

Having chosen a nontrivial epimorphism $\chi$ and associated cover $p:
\tilde M \to M$, we consider  the relative homology group
$H_1(\tilde M, p^{-1}(*))$, which is a finitely generated module over
$\L=\Z\Pi$. We denote the
module by ${\cal A}_t$ and call it the {\it Alexander-Fox
module} of $t$. Generators for ${\cal A}_t$ can be found
corresponding to the arcs of any diagram of $t$; a set of defining
relations is obtained from crossings, as in Figure 2 below. A
similar description for 
oriented links ($0$-tangles) is well known.  Details can be found in
[{\bf SW00}], for example.\bs

\epsfxsize=.8truein
\centerline{\epsfbox{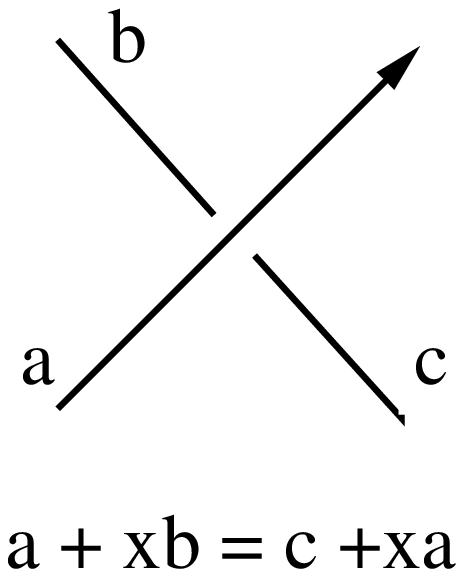}} \bs

\centerline{Figure 2. Crossing relation in ${\cal A}_t$}\bs  

In the case that $M$ is the 2-fold cyclic cover of $B^3$ branched over $t$ and $\chi$ is trivial, we can still obtain $H_1M \oplus {\Bbb Z}$ using the procedure above, letting $x=-1$. Killing the generator corresponding to any single arc yields
$H_1M$.

When $d>1$ the Alexander
module ${\cal A}_t$ is easier to present with generators and relators than is
$H_1\tilde M$.  The two modules fit into
the following exact sequence which arises from the long exact homology sequence of the pair $p^{-1}(*) \subset \tilde M$.
$$0\ \to \ H_1\tilde M \ {\buildrel f \over \longrightarrow}\ {\cal A}_t
\ {\buildrel  g \over \longrightarrow}\ \e(\L) \to\ 0\eqno(3.1)$$
Here $\e(\L)$ is the ideal of $\L$ generated by $x_1-1, \ldots, x_d-1$. 
Since $\e(\L)$ is a submodule of $\L$, it is torsion-free. It follows from exactness that  $TH_1\tilde M
\cong T{\cal A}_t$.  Hence $\t = \D_0(TH_1\tilde M)$ can be computed as
$\D_0(T{\cal A}_t)$. It is not necessary to find the torsion submodule;
Lemma 2.3 (2) ensures that $\t$ is equal to the first nonzero elementary
divisor $\D_k({\cal A}_t)$. 

We can also compute the boundary invariant $\b$ using the Alexander module.
Let $D$ be the submodule of ${\cal A}_t$ generated by its $\L$-torsion elements and by generators associated to input and output arcs of $t$.  By
Lemma 2.3 (4) the sequence
$$0\ \to \ H_1\tilde M/f^{-1}(D) \ {\buildrel \bar f \over \longrightarrow}\ {\cal
A}_t/D
\ {\buildrel \bar g \over \longrightarrow}\ \e(\L)/g(D) \to 0$$
is exact. The preimage
$f^{-1}(D)$ is generated by $TH_1\tilde M$ together with the image of $i_{1+}:H_1 (\p_+ \tilde M)\to H_1\tilde M$. Hence
$H_1\tilde M/f^{-1}(D) \cong BH_1\tilde M/ {\rm im}\ \nabla$ (see Definition 2.1). 
Each generator  is mapped by $g$ to $x_j-1$ in the exact
sequence above, where $j$ corresponds to the component of $t$ to which the
associated arc belongs. Hence $\e(\L)/g(D)$ is either trivial (for example, if $t$ has
no closed components) or else its rank is $1$. In the latter case, we apply
the following. \bs

\ni {\bf Lemma 3.3.} (Lemma 7.2.7(3), [{\bf Kw96}]) Let $0\to A \to B\to C\to 0$ be
a short exact sequence of finitely generated $\L$-modules. If $TC = 0$
and the rank of $C$ is $1$, then $\D_0(A) \= \D_1(B)$. \bs

\ni {\bf Example 3.4.} Consider the colored, oriented ``square-tangle'' $t$
with arcs labeled as in Figure 3 below. The Alexander-Fox module
${\cal A}_t$ has a presentation with generators
$a,b,c,d,e,f,g,h$ and relations: 
$$\eqalign{&b + ya = c + xb,\ c + xb = d + yc,\ d+yc = e + x d, \cr & 
b + hy = g + xb,\  g + xb = f + yg,\ f + yg = e + xf.}$$
(Here we use $x,y$ instead of the more cumbersome $x_1, x_2$.)
By elementary operations we find that $c,d,e,f$ and $g$ can be 
expressed in terms of $a,b$ and $h$. More precisely:
$c= ya + (1-x)b, d = (y-y^2)a + (1-y+xy)b, e = (y-xy+xy^2)a+(1-x+xy-x^2y)b, f =
(1-y+xy)b+(y-y^2)h$ and
$g=(1-x)b+yh$. We have
$${\cal A}_t \cong \<a,b,h \mid (1-x+xy)a = (1-x+xy)h\>.$$
The $\L$-torsion submodule is isomorphic to $\L/(1-x+xy)$.
Hence $\t= 1-x+xy$. We can compute the boundary 
invariant $\b$ by first killing the images of  $a, d,f$ and $h$ in $BH_1\tilde M \cong
\<b\mid\ \>$, obtaining the quotient module $\<b\mid (1-y+xy)b\>$, and then taking
the $0$th elementary divisor. We find that $\b = 1-y+xy$. 

If we reverse the orientation of one component of $t$, say the first, 
then $\t$ and $\b$ become $1-x^{-1}+x^{-1}y \= 1-x-y$ and 
$1-y+x^{-1}y \= x + y -xy$, respectively. \bs

\epsfxsize=1truein
\centerline{\epsfbox{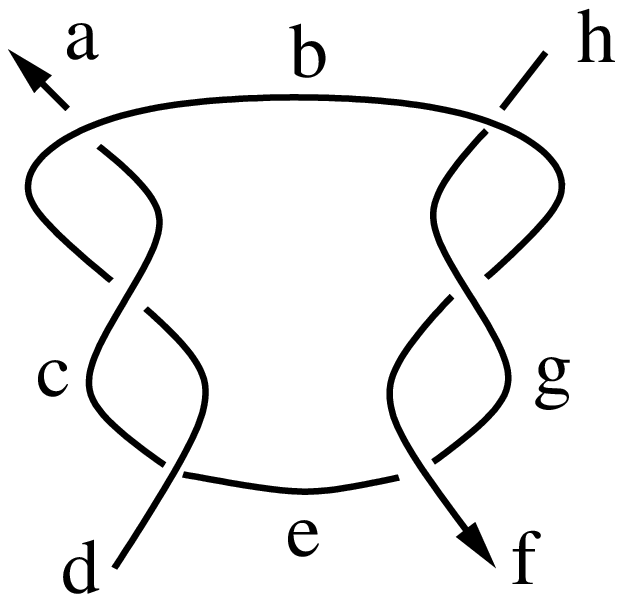}} \bs

\centerline{Figure 3: Labeled tangle $t$}\bs

\ni {\bf Remark 3.5} In the example above $\b$ divides $\bar \t$ (in fact, they are equal). In all examples that we have computed, $\b$  divides $\bar \t$ provided that $\b$ is nonzero.  \bs

Assume that a $2n$-tangle $t$ embeds in a link $\ell=
\ell_1 \cup \cdots \cup \ell_d$. We denote the exterior $S^3 - int\ N(\ell)$ by $E_\ell$. 
If colors and orientations of $t$ (if any) match those of $\ell$, then the augmentation
homomorphism $\chi$ for $t$ extends over $H_1 E_\ell$. We will assume that this is the 
case. Let $\tilde E_\ell$ denote the corresponding cover of $E_\ell$.  

If $t$ and $\ell$ are colored and oriented, and $\chi$ maps the $i$th oriented meridian
to $x_i\in \Pi$, then $\tilde E_\ell$ is the maximal abelian cover of $\ell$; in this case,
$\D_0(H_1\tilde E_\ell)$ is the multivariable Alexander polynomial $\D_l(x_1, \ldots,
x_d)$ of the link.  If $\ell$ is merely oriented, and $\chi$ maps each oriented meridian
to $x \in \Pi\cong \<x\mid \>$, then $\tilde E_\ell$ is an infinite cyclic cover; in this
case, $\D_0(H_1\tilde E_\ell)$ is the $1$-variable Alexander polynomial $\D_l(x)$ of the
link. In either case,  Theorem 2.2 implies that the product $\t\b$ of the torsion
and boundary invariants of $t$ divides the Alexander polynomial
of $\ell$.

The $1$-variable Alexander polynomial of an oriented link 
is related to the multivariable  Alexander polynomial. The following
Lemma is a consequence of Proposition 7.3.10(1) of [{\bf Kw96}].\bs

\ni {\bf Lemma 3.6.} If $\ell$ is an oriented link of $d>1$
components, then $\D_l(x) {\buildrel \cdot \over =} (x-1)
\D_l(x, \cdots, x)$.
\bs

Both the multivariable and single-variable Alexander polynomials can be
found directly from a diagram for $\ell$. Consider the
$\Z\Pi$-module with generators
$a,b,c,\ldots$ corresponding to the arcs of the diagram and
relations associated to the crossings, as in Figure 2. (Here $\Pi$ is 
$\<x_1, \ldots, x_d\mid \>$ or $\<x\mid \>$, depending on which polynomial
is desired.) One builds a
presentation  matrix with columns and rows corresponding 
to generators and relators, respectively. Any submatrix  
obtained by deleting a single row and column is an
{\it Alexander matrix} of $\ell$. Its determinant is
the Alexander polynomial.  Details can be found in [{\bf Li97}], for example. 

Finally, we consider the case that $t$ is neither colored nor oriented.
Let $M$ be the $2$-fold cyclic cover of $B^3$ branched over $t$. When
$t$ embeds in a link $\ell$, then $M$ embeds in $N$, the $2$-fold cyclic
cover of $S^3$ branched over $\ell$. In order to apply Theorem 2.2
we let $\chi$ be the homomorphism mapping $H_1M$ and $H_1N$ to the trivial
group. It is well known that $\D_0(H_1N)$ is the determinant of $\ell$ 
(see for example [{\bf Li97}]). It follows from [{\bf Ru00}] that 
the product $\t\b$ is the greatest common divisor of ${\rm det}\ n(t)$
and ${\rm det}\ d(t)$. Hence Krebes's theorem (Theorem 2.1) is 
a consequence of Theorem 2.2. 
For the convenience of the reader we repeat the argument of 
[{\bf Ru00}]. 

Assume that $t$ is an uncolored, unoriented $4$-tangle that embeds in 
a link $\ell$. Let $M$ be the $2$-fold cyclic cover of $B^3$
branched over $t$. The boundary of $M$ is a torus. Let $N$
be the $2$-fold cyclic cover of $S^3$ branched over $\ell$. Then
$M \subset N$. Moreover, the order of $H_1N$ is the determinant of $\ell$.
Poincar\'e duality implies that
$M$ and $M' = N - {\rm im}\ M$ are rational homology circles. Hence we
can write
$H_1M \cong \Z\oplus \Z/q_1\oplus\cdots\oplus\Z/q_s$, for 
some positive integers $q_1, \ldots, q_s$. Note that the order
of $TH_1M$ is $|q_1\cdots q_s|$. From the long exact homology 
sequence of the pair $M \subset N$:
$$\cdots \to H_2(N,M)\ {\buildrel \p \over \longrightarrow}\ H_1M
\to H_1N \to \cdots$$
we see that $H_1M/\p H_2(N,M)$ embeds in $H_1N$. By the excision
isomorphism $H_2(N,M)\cong H_2(M',\p M)$ is infinite cyclic, generated
by a relative $2$-cycle $C$. The boundary $\p C$ represents
a class $\g \in H_1M$. By what we have already said, $H_1M/\<\g\>$
embeds in
$H_1N$.  Let $i_{1+}: H_1(\p M) \to H_1M$ be the homomorphism induced
by inclusion. Then $i_{1+}C = (c, c_1, \ldots, c_s)$, for some integers
$c, c_1,
\ldots, c_s$. Note that the order of $H_1M/\<\g\>$ is equal to 
$|c|$ times the order of $TH_1M$. 

It is clear that $H_1M/\<\g\>$ is presented by the matrix
$$\pmatrix{c&c_1&\cdots&c_s\cr 0&q_1&0&\cdots& \cr
\vdots&&&\vdots \cr 0 & \cdots &0&q_s\cr}.$$
The class $\g$ is equal to $m_1 \mu + m_2 \l$, where $\mu$ is the
class of the meridian of $\p M$, $\ell$ is the class of the longitude
and $m_1, m_2$ are integers. Certainly gcd $(m_1, m_2) \cdot 
|q_1\ldots q_s|$  divides
the order of $H_1N$. However, $|m_1\cdot q_1\ldots q_s|$ and 
$|m_2\cdot q_1\ldots q_s|$ are easily seen to be the orders
of $H_1N$ when $\ell$ is the numerator and denominator closures
of $t$. Hence ${\rm gcd}\ (m_1, m_2) \cdot 
|q_1\ldots q_s|$ is the greatest common divisor  of the determinants of
$n(t)$ and $d(t)$, and it divides the determinant of $\ell$.

The relation between this proof, due to Ruberman,  and our approach
is the following. 
The absolute value of the boundary invariant $\b$ is the order of
$BH_1M/{\rm im}\ \nabla$. This quotient module of  $H_1M\cong \Z \oplus
\Z/q_1\oplus \cdots \oplus \Z/q_s$ can
be obtained by killing the torsion
elements and then killing the cosets of $i_{1+}\mu$ and $i_{1+}\l$.
The result is a cyclic group of order ${\rm gcd}\ (m_1, m_2)$. Hence
gcd $(m_1, m_2)= |\b|$. The quantity $|q_1\cdots q_s|$ is the
order of $TH_1M$, which is $|\t|$. \bs

\ni {\bf Example 3.7.} We return to the square tangle $t$ of Example
3.4. Ignoring colors we find that $t$ and $\b$ are both equal to
$x^2-x+1$. Hence $(x^2-x+1)^2$ divides the Alexander polynomial
of any oriented link in which $t$ embeds. 

If we ignore the orientation of $t$ as well, then $\t$ and $\b$
are both equal (up to sign) to 3. Hence $9$ divides the determinant 
of any link $\ell$ in which $t$ embeds. 
 \bs

\ni {\bf Example 3.8.} Consider the uncolored, oriented $6$-tangle $t$
in Figure 4. In the exact sequence (3.1) the module $\e(\L)$ is
cyclic, and thus the sequence splits. As a result $H_1\tilde M$
is isomorphic to the quotient ${\cal A}_t^0$  of the Alexander-Fox
module ${\cal A}_t$ obtained by killing a single meridianal
generator. In Figure 4, five meridianal generators are labeled, one
of them with zero; the remaining generators can be written
in terms of these. A presentation matrix $A$ can be found for
${\cal A}_t^0$. For convenience we choose $A$ to be square ($4\times
4$) with two zero rows. We used the software package {\it Maple}, which enabled us
to find nonsingular $4\times 4$ matrices $U,V$ over the ring $R={\Bbb
Q}\Pi$ such that $UAV$ is a diagonal matrix ${\rm
diag}(1,\t(x),0,0)$, the Smith normal form of $A$, where $\t(x)=
(x^2-4x+1)(x^2-x+1)(x-1)$. Hence
$${\cal A}_t^0\otimes_\Z{\Bbb Q} \cong ({\Bbb Q}\Pi)^2\oplus{\Bbb
Q}\Pi/(\t(x)).$$
The polynomial $\t(x)$ is the torsion invariant of the $6$-tangle.
In order to find the boundary invariant $\b(x)$ we project the
five generators of ${\cal A}_t^0$ corresponding to input/output
arcs of $t$ onto the free part of ${\cal A}_t^0\otimes_\Z{\Bbb Q}$,
using the matrix $V$. The images comprise the rows of a $5\times
2$-matrix, and the greatest common divisor of the the $2\times
2$-minors of this matrix is $\b(x)$. We find that $\b(x)=x^2-x+1$.
By Theorem 3.5 the polynomial
$\b\t=(x^2-x+1)^2(x^2-4x+1)(x-1)$ divides the 
Alexander polynomial of any oriented link $\ell$ in which $t$ embeds.
Ignoring orientations, we find that $108$ divides the
determinant of $\ell$. \bs

\epsfxsize=2truein
\centerline{\epsfbox{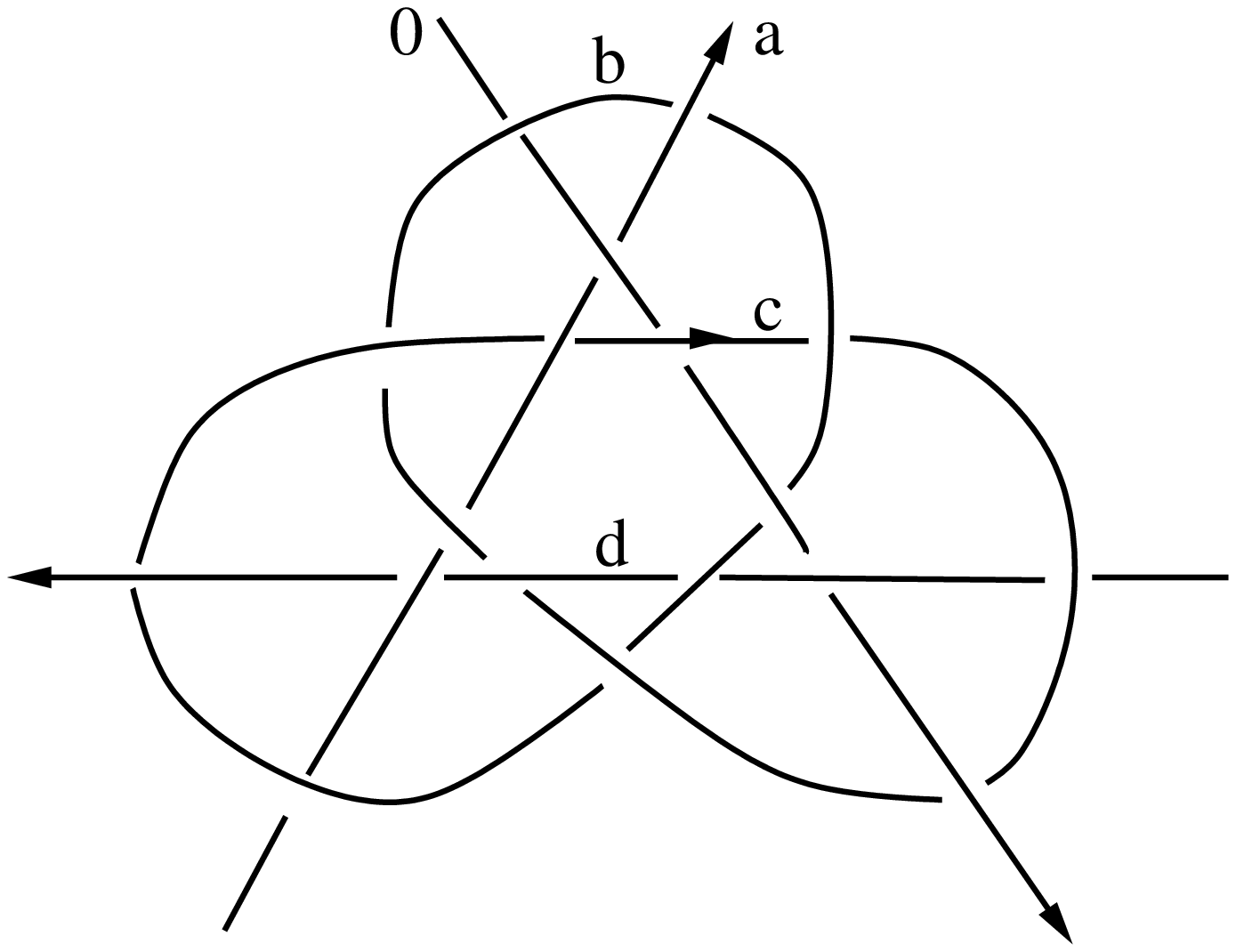}} \bs

\centerline{Figure 4: Labeled $6$-tangle $t$}

\bs


\ni {\bf 4. Virtual links and $2n$-tangles.} In 1996 L.~Kauffman
introduced the notion of a virtual link [{\bf Ka97}], thereby
extending the ``classical'' category of links in an interesting
and nontrivial manner. We review the main ideas. The reader can
find additional information in [{\bf Ka99}] or [{\bf Ka00}]. 

A diagram for a classical link is a a planar $4$-regular graph
with information at each vertex indicating how the link crosses
itself when viewed from a fixed perspective. The decorated vertex
is called a crossing. A {\it virtual link diagram} is likewise
defined. However, such a diagarm is permitted to contain crossings
of a new, ``virtual'' type. Classical and virtual crossing
conventions appear in Figure 5. In many respects virtual crossings
are treated as though they are not present. For example, the arcs
of a diagarm are defined to be the maximal connected components,
just as for classical link diagrams, regardless of the virtual
crossings that they might contain. \bs

\epsfxsize=1.3truein
\centerline{\epsfbox{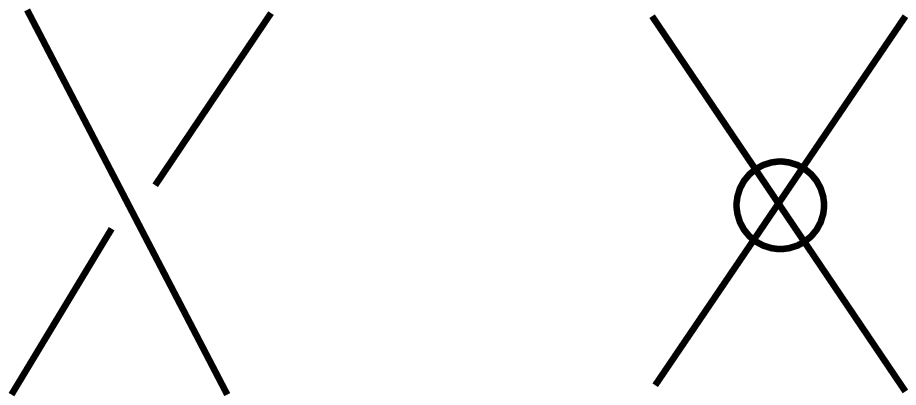}} \bs

\centerline{Figure 5: Classical and virtual crossings} \bs

Two virtual link diagrams are equivalent if one can be obtained
from the other by a finite sequence of the usual, classical
Reidemeister moves or ``virtual'' Reidemeister moves. Virtual
Reidemeister moves are shown in Figure 6. By a generalized
Reidemeister move we mean either a classical or virtual
Reidemeister move.

A {\it virtual link} is an equivalence class of diagrams. We define
{\it virtual $2n$-tangle diagram} and {\it virtual $2n$-tangle} in
the  same manner. Orientations can be imposed as in the classical
category. As before, a virtual tangle means a virtual $4$-tangle.
The numerator closure and the denominator closure are defined as
in the classical case. 

\bs
\epsfxsize=3.4truein
\centerline{\epsfbox{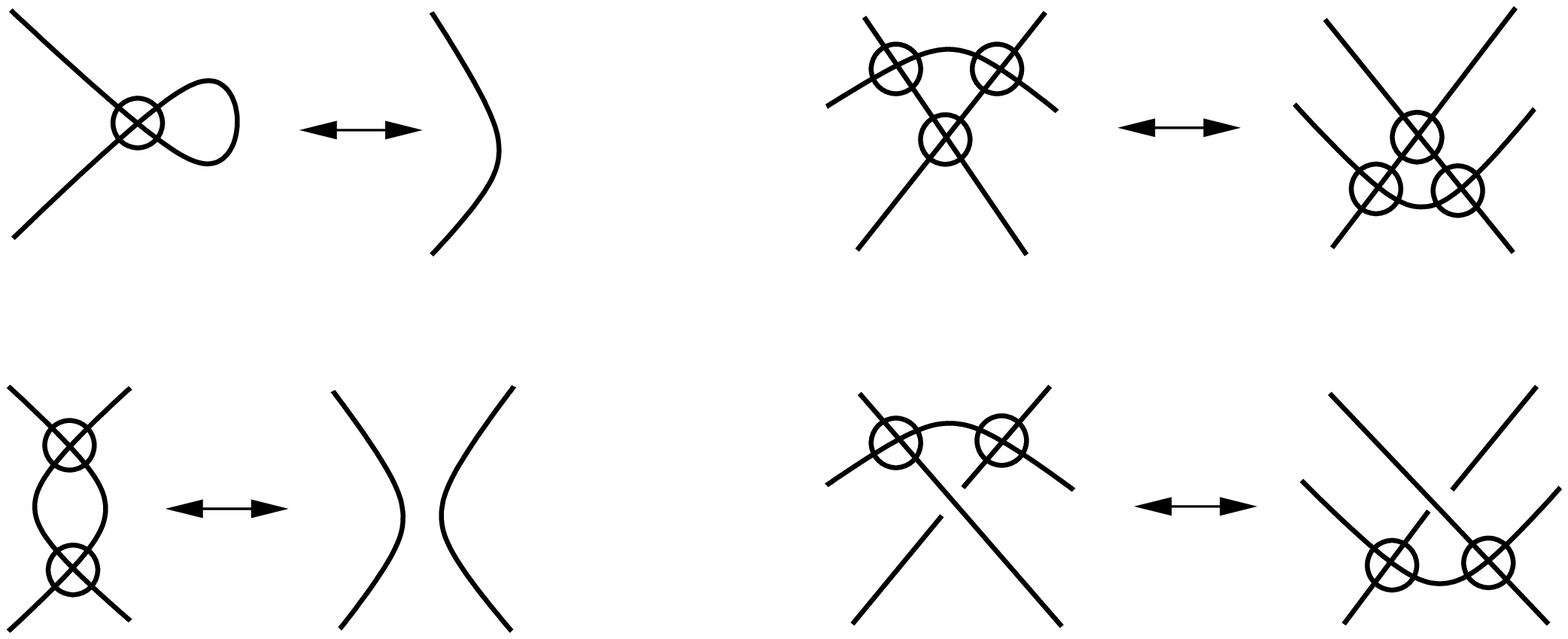}} \bs

\centerline{Figure 6: Virtual Reidemeister moves}\bs

A large body of  virtual knot literature has already appeared. One
reason for the interest is that virtual knot theory contains the
classical theory; more precisely, if two  classical knots are
equivalent under generalized Reidemeister moves, then they are
equivalent under classical ones. A proof can be found in [{\bf
GPV00}] (see also [{\bf Ka}]).

An Alexander matrix can be associated to a diagram of  an oriented
virtual link $\ell$ by the same procedure as in section 3 (see the
paragraph following the proof of Lemma 3.4.) The {\it Alexander
polynomial} $\D_\ell(x)$ is defined to be the greatest common divisor
of the determinants of all submatrices obtained by deleting a single
row and column. 

When $\ell$ is a classical link any row of an
Alexander matrix is a linear combination of the other rows, and
consequently the determinants of any two submatrices differ by a
unit factor. This need no longer be true when $\ell$ is virtual, and
in that case the determinants of {\sl all} the submatrices must be
considered. 

Setting all of the variables equal to $-1$ in the Alexander matrix matrix, and then taking the greatest common divisor of the submatrices produces an integer. Its absolute value is called the {\it determinant} $\det(\ell)$ of the link. It is well known that in the classical case, $\det \ell$ is equal to the absolute value of the Kauffman bracket polynomial of $\ell$ evaluated at a primitive eighth root of unity $\zeta$.
However, for virtual knots and link, such an evaluation might not agree with the determinant; in fact, it need not be an integer. (This anomaly was pointed out in 
[{\bf SW99}].)

A virtual $2n$-tangle $t$ {\it embeds} in a virtual link $\ell$ if
some diagram for $t$ extends to a diagram for $\ell$.
\bs

\ni {\bf Theorem 4.1.} [{\bf SW99}] Assume that $t$ is a virtual
tangle that embeds in a virtual link $\ell$. If $d$ is a square-free
integer that divides the determinants of both the numerator
closure and the denominator closure of $t$, then $d$ divides
$\det(\ell)$. \bs

We give an example to show that the hypothesis that $d$ is
square-free cannot be relaxed. \bs

\ni {\bf Example 4.2.} Figure 7 shows the square-tangle $t$
embedded in a virtual link $\ell$.  Recall that $\det(n(t))=0$ while
$\det(d(t))=9$. An elementary calculation shows that $\det(\ell)=3$. 
Hence the conclusion of Theorem 4.1 does not
hold when $d$ is 9.  However, the evaluation of the Kauffman bracket polynomial of $\ell$ at a primitive eighth root of unit is equal to 9. We will return to this example in Section 5 (see Example 5.6).\bs

\epsfxsize=1.4truein
\centerline{\epsfbox{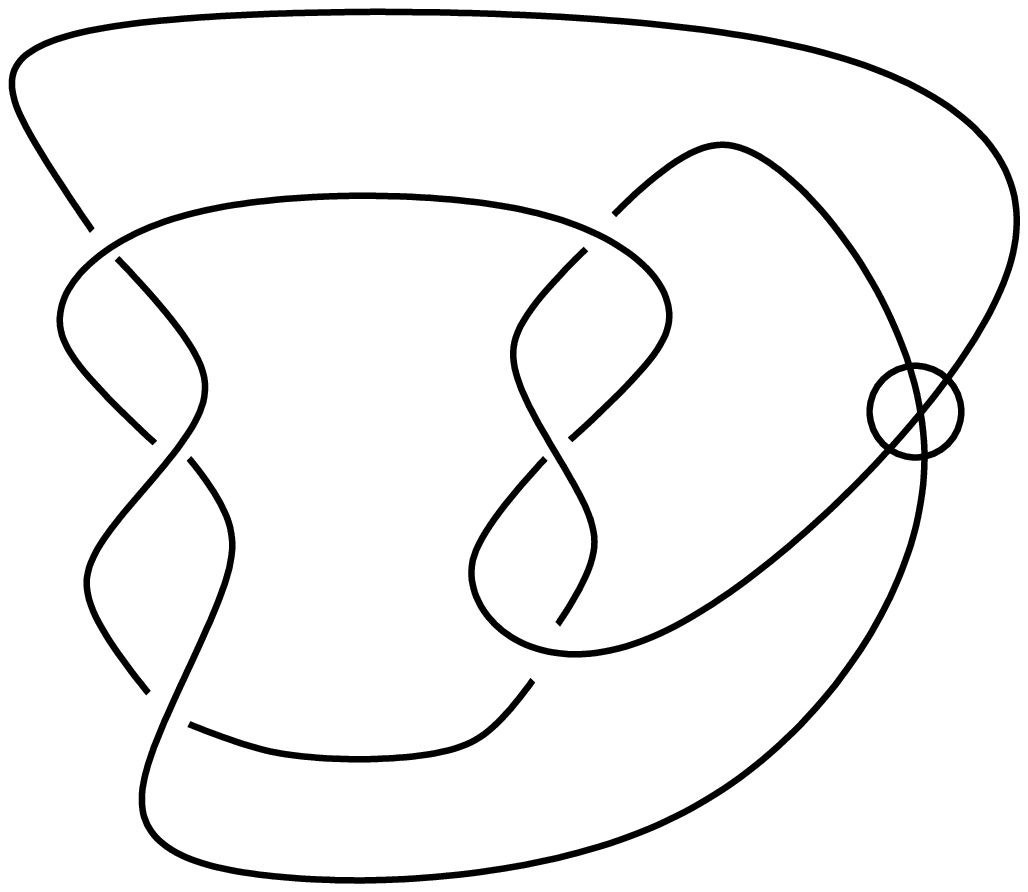}} \bs

\centerline{Figure 7: Embedded square-tangle}\bs

\ni {\bf 5. Persistent invariants from skein theory.} Skein theory and related ideas provide obstructions to embedding tangles. First we discuss applications of the Kauffman bracket skein theory, which allows us to generalize Krebes's theorem to 
$2n$-tangles, for any $n$. Our knot and link notation follows [{\bf Ro76}].

A {\it Catalan tangle} is a $2n$-tangle without crossings or trivial components. There are ${1 \over n+1} {2n \choose n}$ Catalan tangles.  If $t$ and $s$ are $2n$-tangles, then  $t^s$ denotes the link obtained by closing $t$ by $s$; that is, by joining the corresponding ends of $t$ and $s$ without introducing crossings. 
If $t$ is a $4$-tangle and $s$ is the $0$-tangle (respectively $\infty$-tangle) as in Figure 8, then $t^s$ is the numerator closure (respectively, denominator closure) of $t$. Finally $\langle \ell \rangle$ denotes the Kauffman bracket polynomial of a framed link $\ell \subset S^3$. The reader might consult [{\bf Ka91}] for background. \bs

\epsfxsize=2.3truein
\centerline{\epsfbox{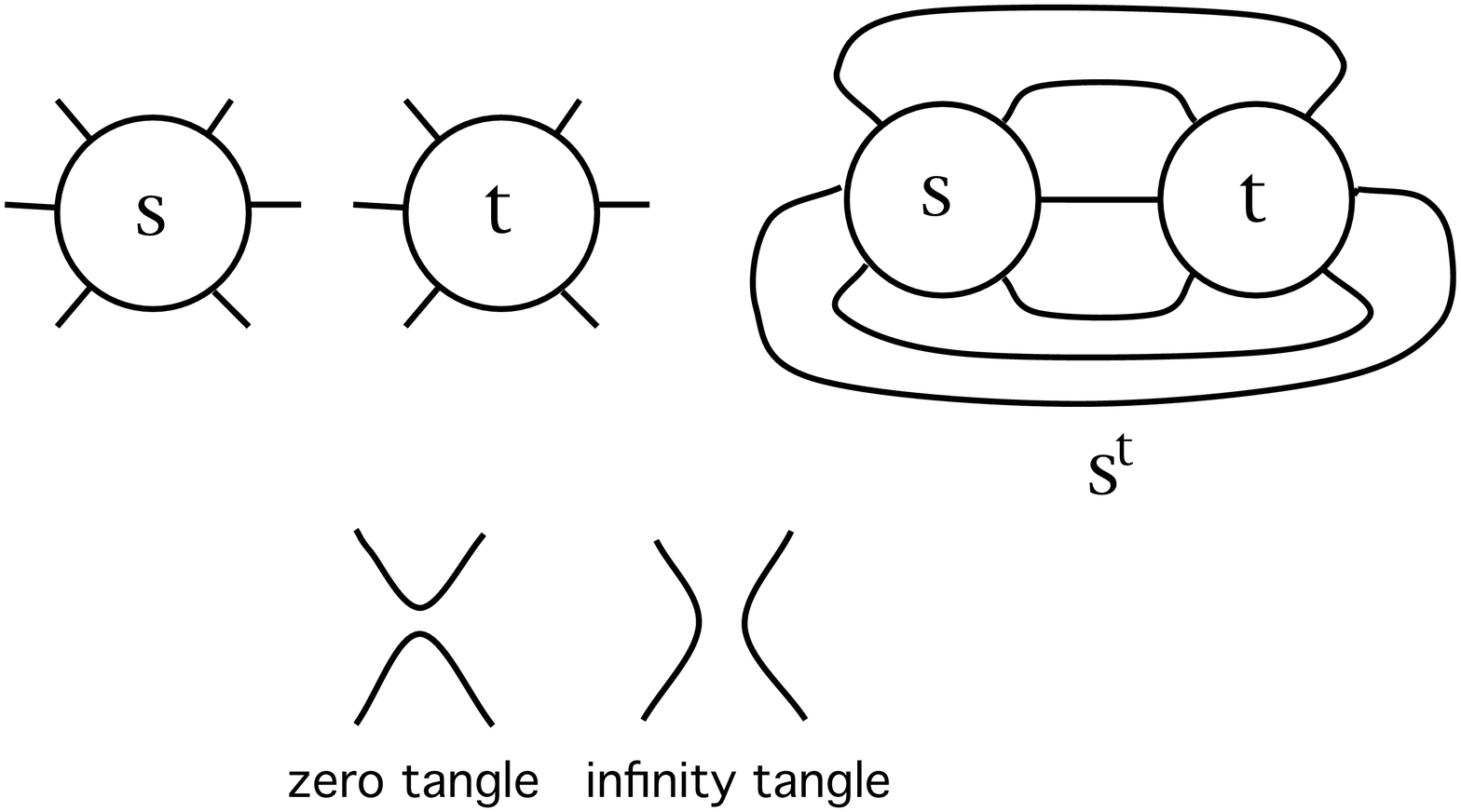}} \bs

\centerline{Figure 8: $2n$-tangle closure $s^t$, $0$-tangle and $\infty$-tangle}\bs

The main result  involving the Kauffman bracket 
polynomial is the following. It offers a technique for deciding whether a $2n$-tangle embeds in a link. 
 \bs

\ni {\bf Theorem 5.1.} If a $2n$-tangle $t$ embeds in a link $\ell$, then the ideal ${\cal{I}}_t$  of ${\Bbb Z}[A^{\pm 1}]$ generated by Kauffman bracket polynomials of all diagrams $\langle t^c\rangle$, where $c$ is any 
Catalan tangle, contains the polynomial  $\langle \ell \rangle$. \bs

\ni {\bf Proof.} Assume that $\ell$ is of the form $t^s$. We use the 
 Kauffman bracket skein relation\footnote{${}^1$}{The relations are $\<\ell_+\> = A\<\ell_0\> + A^{-1}\<\ell_\infty\>$ and $\<\ell \cup \bigcirc \> = (-A^2 -A^{-2})\<\ell\>$.} to eliminate all crossings of $s$, and then eliminate all trivial components of $s$. The resulting 
links are of the form $t^c$, where $c$ are Catalan tangles. As a consequence, 
$\<\ell\>$ is a linear combination of $\<t^c\>$ with the coefficients 
in ${\Bbb Z}[A^{\pm 1}]$. \qed \bs  

We reformulate Theorem 5.1 in the language of Jones polynomials,
recalling that Jones polynomials of the same link with various orientations  
differ only by multiplication by units in
${\Bbb Z}[t^{\pm 1/2}]$. We also use the fact that the determinant 
of the link $\det (\ell)$ is the absolute value of the Jones polynomial evaluated at $t=-1$.\bs

\ni {\bf Corollary 5.2} (i) If a $2n$-tangle $t$ embeds in a link $\ell$, then the ideal
of ${\Bbb Z}[t^{\pm 1/2}]$ generated by Jones polynomials $V_{t^c}$, where $c$ is any
Catalan tangle, contains $V_\ell(t)$.\ss
\ni (ii) The greatest common divisor of the determinants of $t^c$, where $t$ ranges  over all Catalan closures $c$, divides the determinant of $\ell$.  \bs

We illustrate the usefulness of Theorem 5.1 by analyzing the tangle $t$ in
Figure 9, and considering possible links in which it embeds. \bs

\ni {\bf Example 5.3.}  The tangle $t$ in Figure 9 appears in [{\bf Kr99}]. \bs

\epsfxsize=3truein
\centerline{\epsfbox{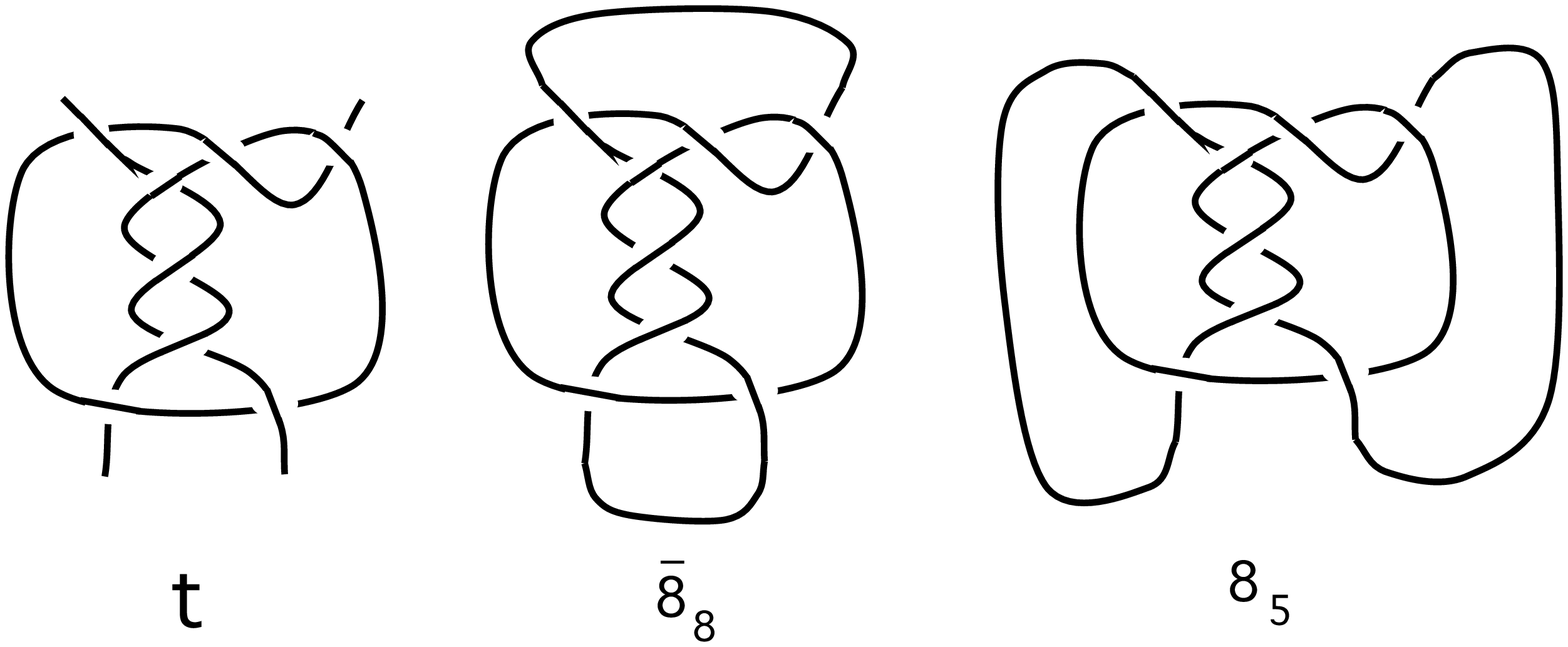}} \bs

\centerline{Figure 9: Krebes's tangle and closures}\bs

\ni The ideal ${\cal I}_t$ is  $(\<\bar 8_8\>, \<8_5\>) = (17, A^4-5)$. 
It is a proper ideal that does not contain  $(A^2+A^{-2})^n$, for any positive integer $n$, nor does it contain $\<4^2_1\>= -A^{10}+A^6-A^2 - A^{-6}$. 
Hence $t$ does not embed in the Hopf link, any trivial link or the link $4^2_1$.  
Furthermore, one can check that for 
knots up to $8$ crossings,  the polynomial $\< k \>$ is contained in the ideal $(17,A^4-5)$ only 
when $k$ is $6_2$, $\bar 8_1$, $\bar 8_{14}$, or of course $8_5$ and $\bar8_8$.
In order to exclude $6_2$ and $\bar 8_1$,  we use 
the criterion based on the Homflypt polynomial (Theorem 5.7). Similarly, 
one excludes $\bar 8_{14}$ with central strands oriented in the same direction.

To find the ideal ${\cal I}_t$ there is no need to consider Catalan 
tangles as in Theorem 5.1.  We used Catalan tangles because they 
form a natural basis of the Kauffman bracket skein module ${\cal S}_{2,\infty}$
of a $2n$-tangle [{\bf Pr91}],[{\bf Pr99}]. Instead we can use any family of $2n$-tangles that generate the skein module, which often allows us to shorten the computation significantly.
In the case of Krebes's tangle $t$, we can replace the numerator $n(t)$ with 
the tangle  $s=$ 
\epsfxsize=.2truein
\epsfbox{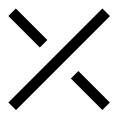}.
Then $t^s$ is the $(4,2)$-torus link $\bar 4^2_1$, with Kauffman bracket polynomial 
$-A^{-10} +A^{-6}-A^{-2} -A^6$. This polynomial is simpler than that of $n(t)$.  We have ${\cal I}_t= (\<\bar 4^2_1\>, \<8_5\>)= (-A^{-10} +A^{-6}-A^{-2} -A^6,\ 
A^{12}-A^8+3A^4-3+3A^{-4}-4A^{-8}+3A^{-12}-2A^{-16}+A^{-20})= (17, A^4-5)$.\bs

Theorem 5.1 generalizes in several directions. For example, we can ask when one tangle embeds in another tangle, as illustrated by Theorem 5.4. The proof is similar to the proof of Theorem 5.1. The product $s\cdot t$ of $2n$-tangles is defined as usual by placing the diagram for $t$ to the right side of a diagram for $s$, and then connecting arcs, as in Figure 10. \bs

\epsfxsize=3truein
\centerline{\epsfbox{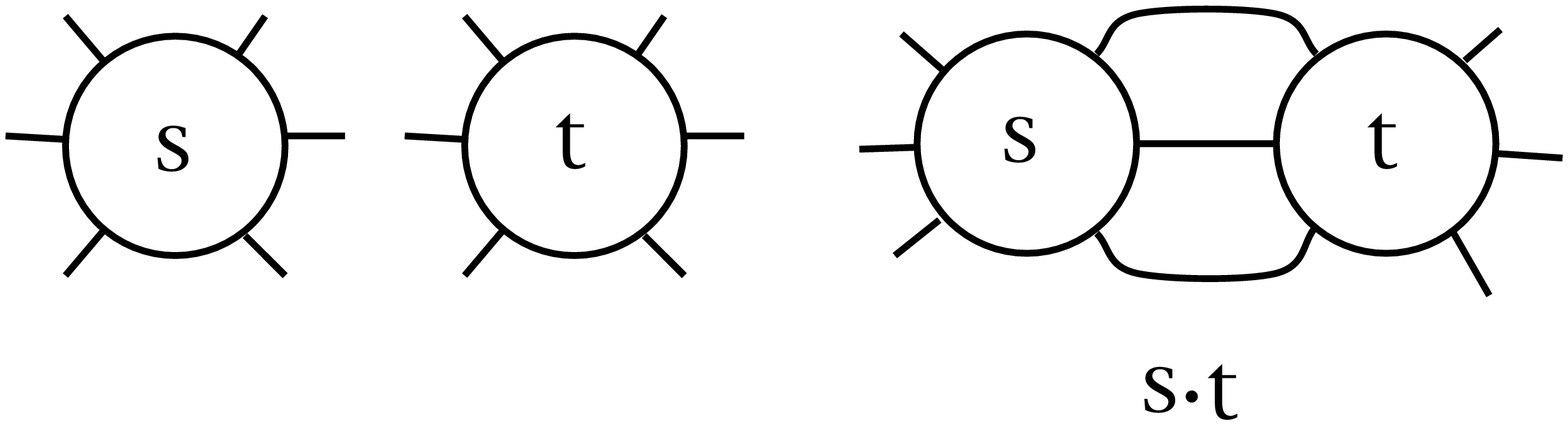}} \bs

\centerline{Figure 10: $2n$-tangle product}\bs

 \ni {\bf Theorem 5.4.} Let $s$ and $t$ be $2n$-tangles. If there exists a $2n$-tangle $u$ such that $t$ is the composition $s\cdot u$, then the element of the Kauffman bracket skein module\footnote{${}^2$}{The Kauffman bracket skein module of a $3$-ball with $2n$ points on its boundary is the quotient of the free module generated by framed unoriented $2n$-tangles modulo the submodule generated by elements $\ell_+ - A \ell_0  - A^{-1} \ell_\infty$ and $\bigcirc \cup \ell + (A^2+A^{-2})\ell.$} represented by $s\cdot u$ is contained in the submodule generated by the elements $s\cdot c$, where $c$ ranges over all Catalan $2n$-tangles. \bs

Another direction extends Theorem 5.1 to the virtual category.  When $n>2$, the  ideal associated to a $2n$-tangle may be larger than in the classical 
case,  as we will illustrate later. The proof of Theorem 5.5 is similar to the proof of Theorem 5.1: we simplify as much as possible the complement $2n$-tangle $t'= \ell - t$ using Kauffman bracket relations extended to the virtual category.  \bs 
 
\ni {\bf Theorem 5.5.}  Let $t$ be a virtual $2n$-tangle and $\ell$ a virtual link.
Let ${\cal I}_t^v$ be the ideal of ${\Bbb Z}[A^{\pm 1}]$ generated by Kauffman bracket polynomials of closures $t^v$, where $v$ ranges over the $(2n)!/2^n n!$  virtual $2n$-tangles corresponding the various ways that the $2n$ boundary points can be connected by $n$ arcs without classical crossings. If $t$ embeds in $\ell$, then $\<\ell\>$ is contained in ${\cal I}_t^v$.\bs

\ni {\bf Example 5.6.} Consider the tangle $t$ in Figure 3. 
The ideal ${\cal I}_t=(A^4+1, (A^{12}+A^4-1)(A^{12}-A^8-1))$ is equal to $(A^4+1, 9)$. The bracket polynomial of the virtual closure $t^v$ is  $A^{20}+A^{18}-A^{16}-2 A^{14}+3A^{10}-2A^6-A^4 +A^2 +1$. Since this polynomial is contained in ${\cal I}_t$, the ideals ${\cal I}_t$ and ${\cal I}_t^v$ are the same in this case.  The single nonclassical closure $t^v$ appears in Figure 7.

Setting $A$ equal to a primitive eighth root of unity $\zeta$, reduces ${\cal I}_t^v$ to the ideal $(9)$ of the ring ${\Bbb Z}[\zeta]$. Hence for any virtual link $\ell$ in which $t$ embeds, the absolute value of $\<\ell\>$ evaluated at $\zeta$ must be divisible by $9$ in ${\Bbb Z}[\zeta]$. (Compare with Example 4.2. See comments preceding Theorem 4.1.) 

\bs

Example 5.6 illustrates a general result: \bs

\ni {\bf Proposition 5.7.} For any classical $4$-tangle $t$, the ideals
${\cal I}_t$ and ${\cal I}_t^v$ are equal. \bs

\ni {\bf Proof.} If a link diagram $D$ has exactly one virtual crossing, denoted 
by $v$, then 
$$(d+1)\<D_v\> = \<D_0\> + \<D_\infty\>,$$
where $d = -A^2 - A^{-2}$. In order to see this, we use Kauffman bracket relations to eliminate all classical crossings and all trivial components not involving $v$. Since $v$ is the only virtual crossing, it suffices to consider only diagrams that are numerator and denominator of the $4$-tangle composed of a single virtual crossing. For such $D$, the skein relation clearly holds. 

We complete the proof of the proposition by applying the above observation to the nonclassical closure $t^v$. We have 
$(d+1)\<t^v\> \in {\cal I}_t = (\<n(t)\>, \<d(t)\>).$
We argue that  $\<t^v\> \in {\cal I}_t$. For this it suffices to show that 
$(d+1, {\cal I}_t)=(1)$ since then $(d+1)\<t^v\> \in {\cal I}_t$ is equivalent to $\<t^v\> \in {\cal I}_t$. 

Since $d+1$ divides $A^8 +A^4+1$, the ideal $(A^8 +A^4+1, {\cal I}_t)$ is contained in $(d+1, {\cal I}_t)$. Notice that ${\cal I}_t$ contains  the bracket polynomial of some classical link $\ell$ with an odd number of components. It follows from statement 12.4 of  [{\bf Jo87}] that $(-A^{-3})^w\<\ell\>-1$ is divisible by $A^8 +A^4 +1$, where $w$ is the sum of the signs of the crossings of a diagram for $\ell$ used to compute $\<\ell\>$. Consequently, $(A^8 +A^4+1, {\cal I}_t) = (1)$ and hence 
$(d+1, {\cal I}_t) = (1)$. 
 \qed\bs

Next we 
generalize Theorem 5.1 by using Homflypt
and Kauffman polynomials in place of the Kauffman bracket polynomial.
The choice of ring is important. 
 
The skein relation of the Homflypt polynomial is 
$$v^{-1}P_{\ell_+}(v,z) -vP_{\ell_-}(v,z) =zP_{\ell_0}(v,z).$$ 
For the initial 
data we take $P_{U_n}=(  {{v^{-1}-v}\over z}  )^{n-1}$, where $U_n$ denotes the trivial link of $n$ components. One might consider
the value of the invariant in the ring ${\Bbb Z}[v^{\pm1},z^{\pm1}]$,
but then every ideal containing $z$ would coincide with the ring. It is
better to use the smaller ring ${\cal R} \subset {\Bbb Z}[v^{\pm 1},z^{\pm 1}]$ 
generated by $v^{\pm 1},z$ and ${v^{-1}-v}\over z$. 
This point of view is used, for example,
in [{\bf Pr89}]. The following result is easily proved. \bs

\ni {\bf Theorem 5.8.} If an oriented $2n$-tangle $t$ embeds in an oriented link $\ell$, then $P_\ell$ is contained in the ideal of ${\cal R}$ generated by polynomials $P_{t^h}$, where $h$ ranges over the $n!$ oriented 
$2n$-tangles that generate the Homflypt skein module. \bs

We have taken the ring ${\cal R}$ instead of
a more familiar ring ${\Bbb Z}[v^{\pm 1},z^{\pm 1}]$ in order to get a stronger
result. For example, the reduction of the ring modulo
$z^k$ leading to Vassiliev invariants is now possible (compare [{\bf Pr94}]).\bs

The question of whether a given element is in an ideal of a polynomial ring can be decided algorithmically using Gr\"obner bases, provided that the coefficient ring is a principal ideal domain.
It can be applied to Krebes's tangle in Example 5.3, for example, 
by using the  Homflypt polynomial and computing the Gr\"obner basis of the associated ideal  in ${\cal R}={\Bbb Z}[v,w,z,y]/(vw-1, zy-v^{-1}+v)$. As another example, the ideal ${\cal J}_t$ for the oriented tangle $t$ in Figure 9 with strands of the central 3-twist oriented in the same direction is generated by $P_{\bar 4^2_1}$ and $P_{8_5}$. From the form of its Gr\"obner basis in ${\cal R}$ we can conclude that $P_{6_2}, P_{\bar 8_1}$ and $P_{8_{14}}$ are not elements of the ideal. We are grateful to M. D{\c a}bkowski for computations in the  ring ${\cal R}\otimes {\Bbb Z}/17{\Bbb Z}$ performed with the program GAP. \bs

\ni {\bf Corollary 5.9.} The Alexander-Conway polynomial $\nabla_\ell(z)$ is contained in the ideal of ${\Bbb Z}[z^{\pm 1}]$ generated by 
elements $\nabla_{t^h}(z)$. \bs

Corollary 5.9 implies that if the $4$-tangle  $t$ can be embedded 
in a link $\ell$, then  $\nabla_\ell$ is in the ideal generated by two polynomials, Alexander-Conway polynomials of the two links obtained from $t$ by closing with its ends using fewest possible crossings. When the orientations of the boundary arcs alternate as one travels along the perimeter of a diagram, these are the numerator and denominator closures; otherwise, one closure acquires an additional crossing while the other does not. 

Corollary 5.9 is useful, as we have seen in Example 3.4. However, it many cases embedding criteria based on the Homflypt polynomial is more helpful. For example, if we try to apply Corollary 5.9 to Krebes's tangle (Figure 9), which we have seen does not embed in many links, we find that the ideal generated by 
$\nabla_{n(t)}$ and $\nabla_{d(t)}$ is trivial, regardless of the orientations chosen. 
 \bs

Similarly we have an obstruction for embedding an unoriented tangle in 
a link using the 2-variable Kauffman polynomial. We use the fact 
that the Kauffman polynomials of links that differ only by orientation of their components are the same up to multiplication by a unit in 
the ring. As before, the ring is taken to be ${\cal R}$.\bs

\ni {\bf Theorem 5.10.} If an unoriented $2n$-tangle $t$ embeds in an unoriented link $\ell$, then $\Lambda_\ell$ is contained in the ideal of ${\cal R}$
 generated by ${\Lambda}_{t^{\kappa}}$, where $\kappa$ ranges over
the $(2n)!\over{2^n n!}$ elements that generate the Kauffman skein module\footnote{${}^3$}
{The Kauffman skein module of a $3$-ball with $2n$ points on its boundary is the quotient of the free module generated by framed unoriented $2n$-tangles modulo the submodule generated by elements $\ell_+  -\ell_-  - \ell_0+\ell_\infty$, $\bigcirc \ \cup\  \ell\  + ({{v^{-1}-v-z}\over z})\ \ell $ and 
$\ell = v\cdot (\ell\ {\rm with\ positive\ twist).}$} 
of $2n$-tangles [{\bf Pr91}]. \bs 

We can look at our criteria from the more general
point of view of bilinear forms on skein modules. Let ${\cal S}(B^3,2n)$
be a skein module of $B^3$ with $2n$ points on its boundary. (See [{\bf HP92}] or [{\bf Pr2}] for details.) It is the quotient of a free module generated by $2n$-tangles by the submodule generated by properly chosen skein expressions.  We have many choices for skein relations.  There is a  bilinear
form $\phi: {\cal S}(B^3,2n) \times {\cal S}(B^3,2n) \to {\cal S}(S^3)$ in which pairs of $2n$-tangles are joined together to form links. The skein module ${\cal S}(S^3)$ is a ring in which product is defined by distant union.
We obtain the following embedding criterion.  If $t$ embeds in $\ell$, then $\ell$ is in the image of the restriction  $\phi(t,\   )$. 
The image is the submodule 
spanned by elements $\phi(t, g)$, where $g$ ranges over a generating set for ${\cal S}(B^3,2n)$.  
For computational
reasons we prefer ${\cal S}(B^3,2n)$ to be finitely generated as is the case for Kauffman bracket, Homflypt and Kauffman skein modules.

\bs


\centerline{\bf References.} 
\bs
\baselineskip=12 pt
\item{[{\bf Bl57}]} R.~C. Blanchfield, {\it Intersection theory of
manifolds  with operators with applications to knot theory}, 
Annals of Math.\ {\bf 65} (1957), 340--356.
\ss
\item{[{\bf GPV00}]} M.~Goussarov, M.~Polyak and O.~Viro,
{\it Finite type invariants of classical and virtual knots},
Topology {\bf 39} (2000), 1045--1068.
\ss

\item{[{\bf HP92}]} J. Hoste and J.H. Przytycki,  A survey of skein modules of $3$-modules, in: Knots 90, Proceedings of the International Conference on Knot Theory and Related Topics, Osaka (A. Kawauchi, ed.), de Gruyter, 1992, 363--379. \ss

\item{[{\bf Jo87}]} V.F.R. Jones, {\it Hecke algebra representations of braid groups and link polynomials},  Annals of Math.\ {\bf 126} (1987), 335--388. \ss

\item{[{\bf Ka91}]} L.H. Kauffman, Knots and Physics, World Scientific, Singapore 1991. 
\ss
\item{[{\bf Ka97}]} L.~H. Kauffman, Talks at MSRI Meeting in
January 1997, AMS Meeting at University of Maryland, College Park
in March 1997. 
\ss

\item{[{\bf Ka99}]} L.~H. Kauffman, {\it An introduction to
virtual knot theory}, European Journal of Combinatorics {\bf 20}
(1999), 663--690.
\ss
\item{[{\bf Ka00}]} L.~H. Kauffman, {\it A Survey of virtual knot theory},
in Knots in Hellas '98, eds. C.McA. Gordon, V.F.R. Jones, L.H. Kauffman, S. Lambropoulou and J.H. Przytycki,World Scientific, Singapore, 2000, 143--202.
\ss
\item{[{\bf Kw96}]} A. Kawauchi, A Survey of Knot Theory,
Birkh\"auser, Basel, 1996. 
\ss
\item{[{\bf Kr99}]} D.~A. Krebes, {\it An obstruction to embedding
4-tangles in links}, J.\ Knot\ Theory\ and\ its\
Ramifications\ {\bf 8} (1999), 321--352. 
\ss
\item{[{\bf KSW00}]} D.~A. Krebes, D.~S. Silver and S.~G. Williams,
{\it Persistent invariants of tangles},  J.\ Knot\ Theory\
and\ its\ Ramifications\ {\bf 9} (2000), 471--477. 
\ss

\item{[{\bf Li97}]} W.~B. Lickorish, An Introduction to Knot
Theory, Springer-Verlag, Berlin, 1997. 
\ss
\item{[{\bf Mi68}]} J.~W. Milnor, {\it Infinite cyclic coverings},
In: Conference on the Topology of Manifolds, J.~G. Hocking, ed., 
Prindle, Weber and Schmidt, 1968. 
\ss
\item{[{\bf Pr89}]} J.H. Przytycki, {\it On Murasugi's and Traczyk's criteria for periodic links}, Math.\ Ann.\ {\bf 283} (1989), 465--478. \ss

\item{[{\bf Pr91}]} J.H. Przytycki, {\it Skein modules of $3$-modules}, Bull.\ 
Acad.\ Pol.\ Math.\ {\bf 39} (1991), 91--100. \ss

\item{[{\bf Pr94}]} J.H. Przytycki, {\it Vassiliev-Gusarov skein modules of $3$-modules and criteria for periodicity of knots}, in: Low-dimensional Topology (Klaus Johannson, ed.) International Press, Cambridge MA 1994, 143--162. \ss

\item{[{\bf Pr99}]} J.H. Przytycki, {\it Fundamentals of Kauffman bracket skein modules}, Kobe\ Math.\ J.\ {\bf 16} (1999), 45--66. \ss

\item{[{\bf Ro76}]} D. Rolfsen, Knots and Links, Publish or
Perish, Berkeley, CA, 1976. 
\ss

\item{[{\bf Ru00}]} D. Ruberman, {\it Embedding tangles in links},
J.\ Knot\ Theory\ and\
its\ Ramifications\ {\bf 9} (2000), 523--530. 
\ss
\item{[{\bf SW99}]} D.~S. Silver and S.~G. Williams, {\it Virtual
tangles and a theorem of Krebes}, J.\ Knot\ Theory\ and\
its\ Ramifications\ {\bf 8} (1999), 941--945. 
\ss

\end